\documentclass[12pt,reqno]{amsart}

\usepackage{amsfonts,amsmath,amsthm}
\usepackage{amssymb,epsfig}
\usepackage{enumerate} 

\usepackage{color} 
\definecolor{vert}{rgb}{0,0.6,0}

\usepackage[text={425pt,650pt},centering]{geometry}
\usepackage{graphicx}
\usepackage{epsfig}
\usepackage{tikz,esvect}
\usetikzlibrary{3d,calc,shapes}

\usepackage{caption}
\usepackage{color} 
\definecolor{vert}{rgb}{0,0.6,0}

\numberwithin{figure}{section}

\theoremstyle{plain}
\newtheorem{thm}{Theorem}[section]

\newtheorem{defn}{Definition}
\newtheorem{quest}{Question}

\newtheorem{lem}[thm]{Lemma}
\newtheorem{cor}[thm]{Corollary}

\theoremstyle{remark}
\newtheorem{rem}{\bf{Remark}}
\numberwithin{equation}{section}

\newcommand{\N}{\mathbb{N}}

\newcommand{\R}{\mathbb{R}}

\newcommand{\T}{\mathbb{T}}
\newcommand{\Z}{\mathbb{Z}}

\newcommand{\AC}{{\rm AC\,}}

\newcommand{\BUC}{{\rm BUC\,}}

\newcommand{\Lip}{{\rm Lip\,}}

\newcommand{\gam}{\gamma}

\newcommand{\ep}{\varepsilon}

\newcommand{\ol}{\overline}

\begin{document}

\title[Differentiability of effective fronts]
{Differentiability of effective fronts in the continuous setting in two dimensions}

\author[H. V. TRAN, Y. YU]
{Hung V. Tran, Yifeng Yu}

\thanks{
The work of HT is partially supported by NSF CAREER grant DMS-1843320 and a Simons Fellowship.
The work of YY is partially supported by NSF grant 2000191.
}

\address[H. V. Tran]
{
Department of Mathematics, 
University of Wisconsin Madison, Van Vleck Hall, 480 Lincoln Drive, Madison, Wisconsin 53706, USA}
\email{hung@math.wisc.edu}

\address[Y. Yu]
{
Department of Mathematics, 
University of California at Irvine, 
California 92697, USA}
\email{yyu1@math.uci.edu}

\date{}

\maketitle

\begin{abstract}
We study the effective front associated with  first-order front propagations in two dimensions ($n=2$) in the periodic setting with continuous coefficients.  
Our main result says that  that  the boundary of  the effective front is differentiable at every irrational point.  
Equivalently,  the stable norm associated with a continuous $\Z^2$-periodic Riemannian metric is differentiable at irrational points.   This conclusion was obtained decades ago  for smooth metrics (\cite{Bangert3, BIK}).     
To the best of our knowledge,  our result provides the first nontrivial property of the effective fronts in the continuous setting, which is the standard assumption  in the PDE theory.  Combining with the sufficiency result in  \cite{JTY}, our result implies that for continuous coefficients, a polygon could be an effective front if and only if it is  centrally symmetric with rational vertices and nonempty interior. 
\end{abstract}


\section{Introduction}
\subsection{Settings}
We first give a minimalistic introduction to the periodic homogenization of Hamilton-Jacobi equations.
For each $\ep>0$, let $u^\ep \in C(\R^n \times [0,\infty))$ be the viscosity solution to 
\begin{equation}\label{C-ep}
\begin{cases}
 u_t^\ep+H\left(\frac{x}{\ep},Du^\ep\right)=0 \qquad &\text{in} \ \R^n \times (0,\infty),\\
u^\ep(x,0)=g(x) \qquad &\text{on} \ \R^n.
\end{cases} 
\end{equation}
Here, the Hamiltonian $H=H(y,p):\R^n \times \R^n \to \R$ is a  given continuous function  satisfying
\begin{equation}\label{con:H}
\begin{cases}
\text{for $p\in \R^n$, $y \mapsto H(y,p)$ is $\Z^n$-periodic;}\\
\text{$H$ is coercive in $p$, that is, } \lim_{|p| \to \infty} \min_{y \in \R^n} H(y,p) = +\infty.
\end{cases}
\end{equation}
For the initial data $g$, we assume
\begin{equation}\label{con:g}
g\in \BUC(\R^n) \cap \Lip(\R^n),
\end{equation}
where  $\BUC(\R^n)$ is the set of bounded, uniformly continuous functions on $\R^n$.

Under assumptions \eqref{con:H}--\eqref{con:g}, $u^\ep$ converges to $u$ locally uniformly on $\R^n \times [0,\infty)$ as $\ep \to 0$, and $u$ solves the effective equation (see \cite{LPV,Ev1,Tran})
\begin{equation}\label{C}
\begin{cases}
 u_t+\ol{H}\left(Du\right)=0 \qquad &\text{in} \ \R^n \times (0,\infty),\\
u(x,0)=g(x) \qquad &\text{on} \ \R^n.
\end{cases} 
\end{equation}

The effective Hamiltonian $\ol{H} \in C(\R^n)$ depends nonlinearly on $H$, and is determined by the cell (ergodic) problems as follows.
\begin{defn}[Effective Hamiltonian]
For each $p \in \R^n$, there exists a unique constant $\ol H(p)\in \R$ such that the following cell problem has a continuous $\Z^n$-periodic viscosity solution
\begin{equation}\label{cell}
H(y,p+Dv)=\ol H(p)  \qquad \text{in} \ \T^n=\R^n/\Z^n.
\end{equation}
Note that $v=v(y,p)$ is not unique even up to additive constants in general.
\end{defn}

If $H$ is convex in $p$, then so is $\ol H$.
In this case, the effective Hamiltonian is also given by an inf-max formula (see, e.g., \cite{Tran})
\begin{align}
\ol H(p)=\inf_{\phi\in C^{\infty}(\Bbb T^n)}\max_{y \in \T^n}\, H(y,p+D\phi(y))
=\inf_{\phi\in C^{1}(\Bbb T^n)}\max_{y \in \T^n}\,H(y, p+D\phi(y)).
\end{align}

It is clear that $\ol H$ is defined in a very implicit way.
A central and fundamental goal in the homogenization theory is to understand qualitative and quantitative properties of $\ol H$.
To date, not much is known about fine properties of $\ol H$.

\smallskip 

We focus on the case $H(y, p)=a(y)|p|$ for $a\in C(\T^n,(0,\infty))$ that arises from the modeling of first-order front propagations  (e.g., crystal growth,  flame propagation), which is probably one of the most physically relevant examples in the homogenization theory.  
In this situation, $\ol H(p)$ represents the effective propagation speed. 
Thanks to the above inf-max formula, 
\begin{align}\label{inf-max}
\ol H(p)=\inf_{\phi\in C^{\infty}(\Bbb T^n)}\max_{y \in \T^n}\,a(y)|p+D\phi(y)|
=\inf_{\phi\in C^{1}(\Bbb T^n)}\max_{y \in \T^n}\,a(y)|p+D\phi(y)|.
\end{align}
Clearly,  $\overline H$ is convex, even, and positively homogeneous of degree $1$. 
We sometime write $\overline H=\overline H_a$ to emphasize the dependence on the function $a$.
Due to those properties of $\overline H_a$, its $1$-sublevel set
\[
S_{a} :=\left\{p\in \R^n \,:\,  \overline H_a(p) \leq 1\right\} 
\]
belongs to $\mathcal{W}$, which denotes the collection of all convex sets in $\R^n$ that are centrally symmetric with nonempty interior.  
The convex dual $D_{a}$ of $S_a$, determined by
\[
D_a=\partial \ol H_a(0),
\]
the subdifferential of  $\ol H_a$ at the origin,
is called the effective front, which also belongs to $\mathcal{W}$.
The following realization problem is of our main interests.

\begin{quest}
For what kind of  $W \in \mathcal{W}$ does there exist a function $a\in C(\T^n,(0,\infty))$ such that  $S_a=W$ (or $D_a=W$)?
 \end{quest}

When $n=2$ and  $a\in C^2(\Bbb T^2, (0, \infty))$,  the following   extra restrictions are known in equivalent forms of stable norms in metric geometry or $\beta$-functions in the Aubry-Mather theory.  See \cite{Bangert1, Bangert3, Carneiro, Mather} for instance. 

\begin{itemize}

\item[(i)] $\partial S_a$ is $C^1$.
Equivalently, $D_a$ is strictly convex.

\item[(ii)] $\partial S_a$ is not strictly convex  (i.e., it  contains   line segments)  unless $a$ is constant. 
 Equivalently,  $\partial D_a$  is not $C^1$  unless $a$ is constant.

\item[(iii)] $\partial S_a$ does not contain a line segment of irrational slope.  
Equivalently, $\partial D_a$ is differentiable at every irrational point. 

\end{itemize}

The $C^2$ regularity was needed in the proofs of (i)--(iii) to ensure that corresponding Hamiltonian systems (or geodesics in the stable norm context) have unique solutions, which implies  that two distinct orbits minimizing the associated actions (or minimal geodesics) cannot intersect twice.  
Together with two dimensional topology, this  leads to a beautiful identification of the minimizing orbits with circle maps that provides a nice characterization of structures of these orbits (see \cite{Bangert1}).
For example, minimal orbits in the same Aubry set are well ordered.

\smallskip

However,  in the merely continuous situation, uniqueness of  solutions to the corresponding ODEs and  intersection restrictions of distinct minimizing orbits cease to exist.  
Two minimizing orbits might intersect multiple (even infinitely many) times. 
Structures of minimizing orbits could be topologically very bad and contain various pathological behaviors.    
As one of the consequences,  it is known now that $\partial S_a$ might not be $C^1$ for $a\in C(\Bbb T^2, (0, \infty)$, i.e., the above property (i) fails.
 In particular,  it was proved in \cite{JTY} that for any $\alpha\in (0,1)$, every polygon with rational slopes in $\mathcal{W}$ can be $S_a$  for some $a$ in $C^{1, \alpha}(\Bbb T^2, (0, \infty))$  (i.e., realizable), which implies that realizable sets are at least dense in $\mathcal{W}$.   It is then tempting to think that every shape in $\mathcal{W}$ might be realizable in the class of $C(\T^2, (0, \infty))$.  
 In this paper, we give a negative answer to this by showing that the above property (iii) still holds in the continuous setting.

\subsection{Main results}
\begin{thm}\label{thm:main-1}
Assume that $n=2$, and $H(y,p)=a(y)|p|$ for $(y,p)\in \T^2 \times \R^2$ for some $a\in C(\T^2, (0, \infty))$. 
Then, $\partial S_a$ does not contain a line segment of irrational slope.  
Equivalently, $\partial D_a$ is differentiable at every irrational point.
\end{thm}

Combining the existence result in \cite{JTY},  this implies that when $n=2$,  for continuous coefficients, a polygon could be an effective front $D_a$ if and only if it is  centrally symmetric with rational vertices and nonempty interior. 

The above result can be proved either from the geometric point of view using minimizing geodesics or  from  the PDE point of view using characteristics of solutions of the cell problem \eqref{cell}.  
In this paper, we choose to use the latter approach.   
For that, let us consider the closely related mechanical Hamiltonian
\[
H(y,p)=\frac{1}{2}|p|^2 + V(y) \quad \text{ for } (y,p)\in \T^n \times \R^n,
\]
for some $V\in C(\T^n)$. Let $\ol H(p)$ be the associated effective Hamiltonian. Note that for $\ol H(p)>\max_{\Bbb T^n}V$, 
\[
\frac{1}{2}|p+Dv|^2 + V(y)=\overline H(p) \quad \Rightarrow \quad  {1\over \sqrt{2(\overline H(p)-V(y))}}|p+Dv|=1.
\]
Accordingly, for $c>\max_{\Bbb T^n}V$ and $a = {1\over \sqrt{2(c-V)}}$, 
\[
F_c =\left\{p\in \R^2\,:\, \ol H(p)=c\right\}=\partial S_a.
\]

Here is our second main result, which implies directly Theorem \ref{thm:main-1}.

\begin{thm}\label{thm:main} 
Assume that $H(y,p)=\frac{1}{2}|p|^2 + V(y)$ for $(y,p)\in \T^n \times \R^n$ for some $V\in C(\T^n)$.
Then, the following properties hold.
\begin{itemize}
\item[(1)] For any $n\in \N$, $p_0, p_1\in \R^n$, if $\overline H(p_0)\not=\overline H(p_1)$, then for all $\lambda\in (0,1)$, 
\[
\ol H(\lambda p_1+(1-\lambda)p_0)<\lambda \ol H(p_0)+(1-\lambda)\ol H(p_1),
\]
that is,  $\ol H$ is strictly convex along directions that are not tangential to the level set in any dimension.

\item[(2)]  For $n=2$ and $c>\max_{\T^2} V$, $F_c$ does not contain a line segment of irrational slope.
\end{itemize}
\end{thm}

Property (1) in the above theorem was  known for smooth $V$ (see \cite {EG}).  
As an immediate corollary, we obtain that 
\begin{cor}\label{cor:lag}
Assume the settings in Theorem \ref{thm:main}. 
Then, the effective Lagrangian $\ol L(q)$, defined as
\[
\ol L(q)=\sup_{p\in \R^2}\left\{p\cdot q-\ol H(p)\right\},
\]
is differentiable at every irrational point.
\end{cor}

\subsection{Connection with stable norms} 
For $a\in C(\Bbb T^n, (0,\infty))$,  we define a corresponding periodic Riemannian metric on $\R^n$ as
\[
g =\frac{1}{a(x)}\sum_{i=1}^{n}dx_{i}^2.
\]  
Let $d_a(\cdot,\cdot)$ denote the distance function induced by this metric.
The \emph{stable norm} associated with $g$ (or $a$) is defined as
\begin{equation}\label{x-a}
\|x\|_{a}=\lim_{\lambda\to \infty}{d_a(0,\lambda x)\over \lambda} \quad \text{ for } x\in \R^n.
\end{equation}
Stable norms and their properties in various settings have been extensively studied in the community of dynamical systems and geometry.  
See \cite{Burago} for more background. 
The limit \eqref{x-a}  can also be viewed as a homogenization problem of the following static Hamilton-Jacobi equation
\[
\begin{cases}
a\left({x\over \ep}\right)|Dw^{\ep}|=1 \quad \text{ for $x\in \R^n\backslash\{0\}$,}\\
w^{\ep}(0)=0.
\end{cases}
\]
Here, $w^\ep$ is the maximal viscosity solution to the above.
By the optimal control formula, for $x\in \R^n$,
\[
w^\ep(x)=\ep d_a\left(0,\frac{x}{\ep}\right).
\]
As $\ep \to 0$, $w^\ep \to w =\|\cdot\|_a$ locally uniformly on $\R^n$, and $w=\|\cdot\|_a$ is the maximal viscosity solution of 
\[
\begin{cases}
\ol H_a(Dw)=1 \quad \text{ for $x\in \R^n\backslash\{0\}$,}\\
w(0)=0.
\end{cases}
\]
The effective front $D_a=\partial \ol H_a(0)$ is then exactly the  unit ball  of  the stable norm $\|\cdot\|_a$
\[
D_a=\{x\in \R^n\,:\,  \|x\|_{a}\leq 1\}.
\]
Our  Theorem \ref{thm:main-1} is equivalent to the conclusion that the stable norm $\|x\|_a=w(x)$ is differentiable at irrational points when $n=2$ for merely continuous metric, which was known for smooth metric \cite{Bangert3}.  See \cite{BIK}  for suitable extensions to higher dimensions.

Somehow,   the word ``homogenization" rarely appeared in the stable norm literature in spite of the aforementioned close connections. 
Similarly, relevant results from stable norms are often not known to people in the PDE community where the focus is more on the wellposedness theory.  
For example,   it was proved in  \cite{Burago} that a periodic Riemannian distance is not much different from the Euclidean distance, i.e., 
\begin{equation}\label{optimal-rate}
|\lambda \|x\|_a-d_a(0,\lambda x)|\leq C \quad \text{ for all } x\in \R^n,
\end{equation}
where $C>0$ is a universal  constant.  
In the PDE language,  by the optimal control formulations, this  is actually equivalent to  the optimal convergence rate $\|w^{\ep}-w\|_{L^\infty(\R^n)}\leq C\ep$ in the static case, and it can be naturally extended to the time-dependent case for general convex Hamilton-Jacobi equations via standard approaches.   
Interestingly, this important fact \eqref{optimal-rate} and its elegant proof  have been hidden from the PDE community for decades until it was found very recently by the authors \cite{TY}.

\subsection{Open problems}  
We list here several open questions that are of interests.

\begin{quest}
 Does there exist a nonconstant $a\in C(\T^2, (0,\infty))$ such that $S_a$ is a strictly convex set (e.g., a disk)? 
 \end{quest}

\begin{quest}
Is any set in $\mathcal{W}$ realizable if we look at $a\in L^{\infty}(\Bbb T^2,(0,\infty))$ with positive essential lower bound?
\end{quest}

\begin{quest}
What can we say about the effective front in higher dimensions ($n\geq 3$), which is much harder due to the lack of topological restrictions?   
A counterexample constructed in  \cite{BIK} says that Theorem \ref{thm:main-1} is not always true when $n\geq 3$ even for smooth $a(y)$. 
In addition, by a clever approach that does not rely on two dimensional topology,  a suitable generalization of Theorem \ref{thm:main-1} was proved in \cite{BIK} for $a\in C^3(\T^n, (0,\infty))$.  
An interesting question is whether the $C^3$ assumption there can be relaxed, which is closely related to regularity of solutions to Hamilton-Jacobi equations. 
See Remark \ref{regularity-3d}. 
In fact, the so called  ``generalized coordinates" in \cite{BIK} are equivalent to viscosity subsolutions $u=p\cdot y +v(y)$ to the cell problem 
\[
a(y)|Du(y)|=a(y)|p+Dv(y)|= 1 \qquad \text{ in } \T^n.
\] 
\end{quest}

\subsection*{Organization of the paper}
The paper is organized as follows.
Some preliminary results are given in Section \ref{sec:prelim}.
We then give the proof of Theorem \ref{thm:main} in Section \ref{sec:proof}.  One of the key ideas is  to rearrange pieces of intersecting minimizing orbits  (Remark \ref{glue}  and Remark \ref{rem:two-curves}) in order to simplify the topology and reduce to situations that are kind of similar to the classical scenarios in \cite{Bangert3}.


\section{Preliminaries} \label{sec:prelim}

In this section, we review some relevant concepts and facts from weak KAM theory \cite{Fa} and Aubry-Mather theory \cite{Bangert1}  in the continuous setting.  
See \cite{FS} for extensions to more general quasiconvex Hamiltonians.  
We assume throughout this section
\[
H(y,p) = \frac{1}{2}|p|^2 + V(y) \quad \text{ for all } (y,p) \in \T^n \times \R^n,
\]
where $V\in C(\T^n)$. 
Then, the cell problem reads
\begin{equation}\label{E-p}
{1\over 2}|p+Dv(y)|^2+V(y)=\ol H(p).
\end{equation}
It is important to note that we only have $V$ is merely continuous on $\T^n$. 	
Except the differentiability property Lemma \ref{regularity},  most of the definitions and proofs are straightforward extensions of those in the smooth setting. 

\smallskip

\begin{defn}\label{ab-curve}
Let $\eta\in \AC([a,b], \R^n)$ be a given curve.
Then, $\eta$ is called an absolute minimizer of the action $\int ({1\over 2}|\dot \gamma(t)|^2-V(\gamma(t))+c)\,dt$ if 
\[
\int_{t_1}^{t_2}\left({1\over 2}|\dot \eta(t)|^2-V(\eta(t))+c\right)\,dt\leq \int_{s_1}^{s_2}\left({1\over 2}|\dot \gamma(t)|^2-V(\gamma(t))+c\right)\,dt
\]
for any $[t_1,t_2]\subset [a,b]$, and $\gamma\in \AC([s_1,s_2], \R^n)$ satisfying that $\gamma(s_i)=\eta(t_i)$ for $i=1,2$.  
Here $\AC(I, \R^n)$ represents the set of absolutely continuous curves defined on the interval $I$. 
\end{defn}

The following lemma is quite basic and well-known (see \cite{Fa,Tran}).

\begin{lem}\label{calib} 
Let $U$ be an open subset of $\R^n$. 
Assume that for some $c\in \R$,  $w\in W^{1,\infty}(U)$ satisfies that 
\[
{1\over 2}|Dw|^2+V(x)\leq c \quad \text{ for a.e. $x\in U$.}
\]
Then, for any  $\eta\in \AC([a,b], U)$,
\[
I[\eta, (a,b)]=\int_{a}^{b} \left({1\over 2}|\dot\eta(t)|^2-V(\eta(t))-c\right)\,dt\geq w(\eta(b))-w(\eta(a)).
\]
\end{lem}

Next is a differentiability property, an important and new result in the merely continuous setting.

\begin{lem}\label{regularity} 
Let $U$ be an open subset of $\R^n$. 
Assume that for some $c\in \R$, $w_i\in W^{1,\infty}(U)$ for $i=1,2$ satisfies that 
\[
{1\over 2}|Dw_i|^2+V(x)\leq c \quad \text{ for a.e. $x\in U$.}
\]
Assume further that there exists a curve $\eta\in \AC([a,b], U)$ such that, for $i=1,2$, 
\[
\int_{a}^{b} \left({1\over 2}|\dot\eta(t)|^2-V(\eta(t))+c\right)\,dt= w_i(\eta(b))-w_i(\eta(a)).
\]
Then, the following properties hold.
\begin{itemize}
\item[(1)] For a.e. $t\in [a,b]$,
\begin{equation}\label{constant-energy}
{1\over 2}|\dot \eta(t)|^2+V(\eta(t))=c.
\end{equation}

\item[(2)] If $\eta$ is differentiable at $t_0\in (a,b)$, then $w_1$ and $w_2$ are differentiable at $x=\eta(t_0)$ and 
\[
Dw_1(\eta(t_0))=Dw_2(\eta(t_0))=\dot \eta(t_0).
\]

\item[(3)]    For all $x\in \eta((a,b))$,   $w_1-w_2$ is differentiable  at $x$  and 
\[
D(w_1-w_2)(x)=0.
\]
\end{itemize}
\end{lem}

\begin{proof}
Without loss of generality, we assume that $0\in (a,b)$ and $\eta(0)=0$.

\smallskip

We first prove \eqref{constant-energy}. 
In fact, by standard mollification of $w_1$ and approximations,  we have that 
\[
w_1(\eta(b))-w_1(\eta(a))=\int_{a}^{b}p_1(t)\cdot \dot \eta(t)\,dt=\int_{a}^{b} {1\over 2}\left( |p_1(t)|^2+|\dot \eta(t)|^2-|p_1(t)-\dot \eta(t)|^2\right)\,dt
\]
for some $p_1(t)\in \partial w_1(\eta(t))$ for $t\in (a,b)$.  
Here,
\[
\partial w_1(x)=\mathrm{co}(K(x)),
\]
where ${\rm co}(K(x))$ is the convex hull of the set 
\[
K(x)=\left\{p\in \R^n\,:\, \text{$ \exists \, \{x_k\} \to x$ s.t. $Dw_1(x_k)$ exists, and $p=\lim_{k \to \infty}Dw_1(x_k)$}\right\}.
\]
Apparently, 
\[
{1\over 2}|p_1(t)|^2+V(\eta(t))\leq c  \quad \text{for all $t\in  [a,b]$}. 
\]
Accordingly,  we must have that
\[
{1\over 2}|p_1(t)|^2+V(\eta(t))= c \quad \text{ and } \quad p_1(t)=\dot \eta(t)  \quad \text{for a.e. $t\in  [a,b]$}. 
\]
Hence \eqref{constant-energy} holds.  
Note that this also implies that $\eta$ is Lipschitz continuous. 

\medskip

To prove (2) and (3),  it suffices to show that  for a sequence $\{\lambda_m\} \subset (0,\infty)$ with $\lim_{m\to \infty}\lambda_m=0$, if 
\[
\lim_{m\to \infty}{\eta(\lambda_m t)\over \lambda_m}  \quad \text{ exists for all $t\in \R$},
\]
then,   for all $t\in \R$, $x\in \R^n$, and $i=1,2$,
\[
\lim_{m\to \infty}{\eta(\lambda_m t)\over \lambda_m}=qt \quad \text{ and } \quad \lim_{m\to \infty} {w_i(\lambda_mx)-w_i(0)\over \lambda_m}=q\cdot x
\]
for some $q$ satisfying ${1\over 2}|q|^2+V(0)=c$. 
 We only need to  prove this claim for $w_1$ as the proof for $w_2$ is the same.  
 Let 
\[
\bar \eta(t)=\lim_{m\to \infty}{\eta(\lambda_m t)\over \lambda_m}.
\]
Owing to Lemma \ref{calib},  we have that, for any $a\leq t_1<t_2\leq b$, 
\[
w_1(\eta(t_2))-w_1(\eta(t_1))=\int_{t_1}^{t_2}\left({1\over 2}|\dot \eta(t)|^2-V(\eta(t))+c\right)\,dt.
\]
In fact, $\eta$ is an absolute minimizer of the action $\int ({1\over 2}|\dot \gamma(t)|^2-V(\gamma(t))+c)\,dt$ (see Definition \ref{ab-curve}). Then, ${\eta(\lambda_m t)/\lambda_m}$ is an absolute minimizer of the action 
 \[
 \int \left({1\over 2}|\dot \gamma(t)|^2-V(\lambda_m\gamma(t))+c\right)\,dt
 \]
 over  $\left[{a/ \lambda_m},{b/ \lambda_m}\right]$.
By the stability of minimizing curves,  $\bar \eta$ is an absolute minimizer of the action 
\[
\int \left({1\over 2}|\dot \gamma(t)|^2-V(0)+c\right)\,dt
\]
 over any finite interval of $\R$. 
By the Euler-Lagrange equations, $\bar \eta(t)$ must be a line passing through the origin $0$, that is,
\[
\bar \eta (t)=q\cdot t
\]
 for some $q\in \R^n$.  
 For convenience, denote $M=\sqrt{2(c-V(0))}$.   
 Due to \eqref{constant-energy},  $|q|\leq M$. 
 
By passing to a subsequence if necessary, we assume that 
\[
\lim_{m\to \infty} {w_1(\lambda_mx)-w_1(0)\over \lambda_m}=u(x) \quad \text{ for all } x\in \R^n.
\]
Then, $u\in W^{1, \infty}(\R^n)$, and  $u$ satisfies that 
\[
{1\over 2}|Du(x)|^2+V(0)\leq c  \quad \text{ for a.e. $x\in \R^n$}.
\]
Equivalently, $|Du(x)|\leq M$ for  a.e. $x\in \R^n$.  
Thanks to  \eqref{constant-energy},
\begin{align*}
w_1(\eta(\lambda_mt))-w_1(\eta(0))&=\int_{0}^{\lambda_m t}\left({1\over  2}|\dot \eta(s)|^2-V(\eta(s))+c\right)\,ds\\
&=\int_{0}^{\lambda_m t}  2\left(c-V(\eta(s))\right)\,ds.
\end{align*}
Dividing  both sides by $\lambda_m$, and sending $m\to \infty$,  we derive that 
\[
u(qt)=tM^2 \quad \text{ for all $t\in \R$}.
\]
Since $|Du|\leq M$, 
\[
u(qt)\leq M|q|t\leq M^2t.
\]
Therefore, we deduce that  $|q|=M$, and 
\[
u(et)=Mt  \quad \text{ for all $t\in \R$}
\]
for $e={q\over |q|}={q\over M}$.  
Then, the standard  tightness argument in \cite{Aronsson} leads to
\begin{equation}\label{eq:u-lin}
u(x)=q\cdot x \quad \text{ for all  $x\in  \R^n$.}
\end{equation}
Let us give a proof \eqref{eq:u-lin} here for completeness.  
Due to $|Du|\leq M$, 
\[
|u(x)-Mt|=|u(x)-u(et)|\leq M|x-et|  \quad \text{ for all  $x\in \R^n$, $t\in \R$}.
\]
Taking square of both sides leads to
\[
u(x)^2-2Mtu(x)+M^2t^2\leq M^2(|x|^2-2x\cdot et+t^2), 
\]
which is reduced to 
\[
u(x)^2- M^2 |x|^2 \leq 2Mt(u(x)-q\cdot x).
\]
Since for fixed $x\in \R^n$,  the above holds for all $t\in \R$,  we must have 
\[
u(x)=q\cdot x.
\]
\end{proof}

\begin{rem}\label{regularity-3d}
In the setting of the above lemma, for smooth $V$, it is well known that $\eta$ is $C^2$ and both $w_1$ and $w_2$ are differentiable along $\eta$. See \cite{FS} for differentiability results for Lipschitz continuous coefficients.   In our setting, the subtle point is that it is not very clear to us whether $w_1$ and $w_2$ are  differentiable along $\eta$.  However, $w_1-w_2$ is indeed differentiable along $\eta$ with $D(w_1-w_2)=0$, which is enough for our purpose.   It remains an interesting question to prove or disprove the differentiability of $w_i$ along $\eta$. In particular,  if a proper uniform $C^1$ regularity of $s(x)=w_1-w_2$ near $\eta$ (i.e., $|s(x)-s(\eta(0))|=o(d(x,\eta))$) can be established,  the $C^3$ assumption of $a(x)$ in \cite{BIK} might  be relaxed to the mere continuity assumption via suitable adjustment of the methods there.  
\end{rem}

For $t>0$, and $p, x,y\in \R^n$, let
\[
G_{t,p}(x,y)=\inf_{\substack{\xi\in \AC([0,t],\R^n),\\ \xi(0)=x,\  \xi(t)\in y+\Z^n}}\int_{0}^{t}\left({1\over 2}|\dot \gamma(s)|^2-V(\gamma(s))-p\cdot \dot \xi(s)+\overline H(p)\right)\,ds,
\]
and
\[
G_p(x,y)=\liminf_{t\to \infty}G_{t,p}(x,y).
\]
By Lemma \ref{calib},   $G_{t,p}(x,y)\geq v(y)-v(x)$. 
Hence,
\[
G_p(x,y)\geq v(y)-v(x)
\]
for any viscosity solution $v$ of \eqref{E-p}.  
In particular, $G_p(x,x)\geq 0$.  
In addition,  it is easy to see that $G_p(x,y)$ is $\Z^n$-periodic and Lipschitz continuous in $x$ and $y$.

Now we define the Aubry set associated with $p\in \R^n$ as 
\[
\mathcal {A}_p=\{x\in \R^n\,:\,  G_{p}(x,x)=0\}.
\]
\begin{defn}
Given $p\in \R^n$,   a Lipschitz continuous curve $\xi:\R \to \R^n$ is called a global characteristic associated with a viscosity solution $v$ of the cell problem \eqref{E-p} if, for  all $t_1<t_2$,
\[
v(\xi(t_2))-v(\xi(t_1))=\int_{t_1}^{t_2}\left({1\over 2}|\dot \xi(s)|^2-V(\xi(s))-p\cdot \dot \xi(s)+\overline H(p)\right)\,ds, 
\]
or equivalently, for $u(x)=p\cdot x+v(x)$ for $x\in \R^n$, 
\[
u(\xi(t_2))-u(\xi(t_1))=\int_{t_1}^{t_2}\left({1\over 2}|\dot \xi(s)|^2-V(\xi(s))+\overline H(p)\right)\,ds.
\]
\end{defn}

By Lemma \ref{calib},  to show that $\xi$ is a global characteristics of $v$ is sufficient to show that there exist  $T_m\to +\infty$ and $T_{m}^{'}\to -\infty$ as $m\to +\infty$ such that 
\[
v(\xi(T_m))-v(\xi(T_{m}^{'}))=\int_{T_{m}^{'}}^{T_m}\left({1\over 2}|\dot \xi(s)|^2-V(\xi(s))-p\cdot \dot \xi(s)+\overline H(p)\right)\,ds, 
\]

In addition, owing to Lemma \ref{regularity},  if $\xi$ is a global characteristic associated with $v$ solving \eqref{E-p} for a given $p\in \R^n$, then
\begin{equation}\label{contant-energy-1}
{1\over 2}|\dot \xi(s)|^2+V(\xi(s))=\overline H(p)  \quad \text{ for a.e. $t\in \R$.}
\end{equation}

\begin{rem}\label{abs} 
By Lemma \ref{calib},  for $p\in  \R^n$ and $v$ is a viscosity solution of \eqref{E-p},  every global characteristic $\xi$  associated with  $v$ is an absolute minimizer of the action $\int ({1\over 2}|\dot \gamma(t)|^2-V(\gamma(t))+\ol H(p))\,dt$. 
See Definition \ref{ab-curve}.
\end{rem}

\begin{lem}\label{stability} 
For $p\in  \R^n$ and $v$ is a viscosity solution of \eqref{E-p},  suppose that $\{\xi_m\}$ is a sequence of global characteristics associated with $v$ such that 
\[
\lim_{m\to \infty}\xi_m=\xi  \quad \text{ locally uniformly in $\R$.}
\]
Then, $\xi$ is also a global characteristic of $v$. 
\end{lem}

\begin{proof}
Fix  $t_1<t_2$.  
For $u(x)=p\cdot x+v(x)$ for $x\in \R^n$, and $m\in  \N$,  we have that 
\[
u(\xi_m(t_2))-u(\xi_m(t_1))=\int_{t_1}^{t_2}\left({1\over 2}|\dot \xi_m(s)|^2-V(\xi(s))+\overline H(p)\right)\,ds.
\]
Sending $m\to \infty$, by the lower semicontinuity of the integral, we have that 
\begin{align*}
u(\xi(t_2))-u(\xi(t_1))&\geq \liminf_{m\to \infty} \int_{t_1}^{t_2}\left({1\over 2}|\dot \xi_m(s)|^2-V(\xi_m(s))+\overline H(p)\right)\,ds\\
&\geq \int_{t_1}^{t_2}\left({1\over 2}|\dot \xi(s)|^2-V(\xi(s))+\overline H(p)\right)\,ds.
\end{align*}
Combining this with Lemma \ref{calib}, we get the desired result.
\end{proof}

\begin{defn}
Given $p\in \R^n$,   a Lipschitz continuous curve $\xi:\R \to \R^n$ is called a universal global characteristic associated with $p$  if it is a global characteristic  for every  viscosity solution $v$ of the cell problem \eqref{E-p}.

\smallskip

We denote by  
\[
\mathcal {U}_p=\text{the collection of all universal characteristics associated with $p$. }
\]

\end{defn}

Then for every $\xi\in \mathcal{U}_p$,
\begin{equation}\label{constant-energy-aubry}
{1\over 2}|\dot \xi(s)|^2+V(\xi(s))=\overline H(p)  \quad \text{ for a.e. $t\in \R$.}
\end{equation}
Also, owing to Lemma \ref{stability},  the set $\mathcal{U}_p$ is closed for limit of orbits. 

\begin{lem} 
Given $p\in \R^n$,   every viscosity solution $v$ of \eqref{E-p} has a global characteristic. 
\end{lem}

\begin{proof}
Let $u(x)=p\cdot x+v(x)$ for $x\in \R^n$.  
Then, $w(x,t) = u(x) - \ol H(p)t$ for $(x,t) \in \R^n \times [0,\infty)$ solves the following Cauchy problem
\[
\begin{cases}
w_t + \frac{1}{2}|Dw|^2 + V(x) =0 \quad &\text{ in } \R^n \times (0,\infty),\\
w(x,0) = u(x) \quad &\text{ on } \R^n.
\end{cases}
\]
By the optimal control formula, for each $k\in \N$,
\[
w(0,k)=\inf_{\substack{\gam \in\AC([0,k],\R^n),\\ \gam(k)=0}} \left\{\int_{0}^{k}\left({1\over 2}|\dot \gam(s)|^2-V(\gam(s))\right)\,ds + u(\gam(0)) \right\}.
\]
There exists $\xi_k \in \AC([-k,0],\R^n)$ with $\xi_k(0)=0$ such that
\[
w(0,k)=\int_{-k}^{0}\left({1\over 2}|\dot \xi_k(s)|^2-V(\xi_k(s))\right)\,ds + u(\xi_k(-k)).
\]
This is equivalent to 
\[
u(\xi_k(0))-u(\xi_k(-k))=\int_{-k}^{0}\left({1\over 2}|\dot \xi_k(s)|^2-V(\xi_k(s))+\ol H(p)\right)\,ds
\]
Owing to Lemma \ref{regularity}, 
\[
{1\over 2}|\dot \xi_k(s)|^2+V(\xi_k(s))=\ol H(p) \quad \text{for a.e. $t\in [-k,0]$}.
\]
By sending $k \to \infty$ and passing to a subsequence if necessary, we have that $\xi_k \to \xi$ locally uniformly on $(-\infty,0]$, and, similar to Lemma \ref{stability},  $\xi:(-\infty,0] \to \R^n$ with $\xi(0)=0$ is a backward characteristic associated with $v$.
More precisely, for any $t_1<t_2 \leq 0$,
\[
u(\xi(t_2))-u(\xi(t_1))=\int_{t_1}^{t_2}\left({1\over 2}|\dot \xi(s)|^2-V(\xi(s))+\overline H(p)\right)\,ds.
\]

Again, due to Lemma \ref{regularity},  for a.e. $t\leq 0$, 
\[
{1\over 2}|\dot \xi(t)|^2+V(\xi(t))=\ol H(p).
\]

Next, we create a global characteristic from this backward characteristic $\xi$.
For $m\in \N$, let $\gam_m:(-\infty,m] \to \R^n$ be such that
\[
\gam_m(t) = \xi(t-m) + k_m \quad \text{ for all } t \leq m,
\]
where $k_m\in \Z^n$ is chosen so that $\gam_m(0) = \xi(-m)+k_m \in [0,1]^n$.
Again, by passing to a subsequence if necessary, $\gam_m \to \gam$ locally uniformly on $\R$.
We obtain that $\gam$ is a global characteristic associated with $v$.
\end{proof}

\begin{lem}  
For any $p\in \R^n$, 
\[
\mathcal {A}_p\not=\emptyset.
\]
\end{lem}

\begin{proof}
 Let $v$ be a viscosity solution of \eqref{E-p} and $\xi:\R\to \R^n$ be a global characteristic associated with $v$.  
 By projecting $\xi$ to $\T^n$ and a suitable translation in time, we may find  sequences $\{t_{m}\}\to +\infty$  and $\{x_m\} \subset [0,1]^n$ such that  
\[
\lim_{m\to +\infty}(t_{m+1}-t_m)=+\infty,
\]
\[
  \xi(t_m)=x_m+k_m \quad \text{ for some $k_m\in \Z^n$,}
 \]
 and
\[
\lim_{m\to \infty}x_m=x_0\in [0,1]^n,
\]
and
\[
v(x_{m+1})-v(x_m)=\int_{t_{m}}^{t_{m+1}}\left({1\over 2}|\dot \xi(s)|^2-V(\xi(s))-p\cdot \dot \xi(s)+\overline H(p)\right)\,ds.
\]
Then,
\[
G_{t_{m+1}-t_m,p}(x_m,x_{m+1}) = v(x_{m+1})-v(x_m).
\]
This implies
\[
G_p(x_0,x_0) \leq \liminf_{m\to \infty}G_{t_{m+1}-t_m,p}(x_m,y_m) = 0.
\]
Thus, $x_0\in \mathcal{A}_p$.

\end{proof}

\begin{lem}\label{Aubry-orbit} 
For any $x\in \mathcal{A}_p$, there exists $\xi\in \mathcal{U}_p$ such that $\xi(0)=x$.  
In particular, this implies that $\mathcal{U}_p\not=\emptyset$. 
\end{lem}

\begin{proof}
According to the definition of $ \mathcal{A}_p$, there exist $\{t_m\}\to \infty$ and a sequence of curves $\gamma_m: [0,t_m]\to \R^n$ such that $\gamma_m(0)=x$, and $\gamma_m(t_m)= x+k_m$ for some $k_m \in \Z^n$, and
\[
\lim_{m\to \infty}\int_{0}^{t_m}\left({1\over 2}|\dot \gamma_m(s)|^2-V(\gamma_m(s))-p\cdot \dot \gamma_m(s)+\overline H(p)\right)\,ds=0.
\]
Given a viscosity solution $v$ of \eqref{E-p},  owing to Lemma \ref{calib},  for any fixed $L>0$,  
\[
v(\gamma_m(L))-v(x))\leq \int_{0}^{L}\left({1\over 2}|\dot \gamma_m(s)|^2-V(\gamma_m(s))-p\cdot \dot \gamma_m(s)+\overline H(p)\right)\,ds
\]
and
\[
v(\gamma_m(t_m))-v(\gamma_m(L)))\leq \int_{L}^{t_m}\left({1\over 2}|\dot \gamma_m(s)|^2-V(\gamma_m(s))-p\cdot \dot \gamma_m(s)+\overline H(p)\right)\,ds.
\]
Together with $0=v(\gamma_m(L))-v(x))+v(\gamma_m(t_m))-v(\gamma_m(L)))$,  we have that 
\[
\lim_{m\to \infty}\left(\int_{0}^{L}\left({1\over 2}|\dot \gamma_m(s)|^2-V(\gamma_m(s))-p\cdot \dot \gamma_m(s)+\overline H(p)\right)\,ds-(v(\gamma_m(L))-v(x))\right)=0.
\]
Similarly, 
\begin{multline*}
\lim_{m\to \infty}\Big(\int_{t_m-L}^{t_m}\left({1\over 2}|\dot \gamma_m(s)|^2-V(\gamma_m(s))-p\cdot \dot \gamma_m(s)+\overline H(p)\right)\,ds\\
-(v(\gamma_m(t_m))-v(\gamma(t_m-L)))\Big)=0.
\end{multline*}
Define $\xi_m:\left[-{t_m/2}, {t_m/2}\right]\to \R^n$ as 
\[
\xi_m(t)=
\begin{cases}
\gamma_m(t)  \quad &\text{ for $t\in \left[0, {t_m\over 2}\right]$},\\
\gamma_m\left(t_m+t\right) -k_m  \quad &\text{ for $t\in \left[-{t_m\over 2}, 0\right]$}
\end{cases}
\]
Clearly, for any fixed $L>0$,  $\{\|\xi_m\|_{H^1((-L,L))}\}$ is uniformly bounded.  
Up to a subsequence if necessary, we may assume that 
\[
\lim_{m\to \infty}\xi_m(t)=\xi(t)   \quad \text{ locally uniformly in $\R$.}
\]
for $\xi\in \AC(\R, \R^n)$. 

Then, using Lemma \ref{calib}  and the lower semicontinuity of the integral (similar to the proof of Lemma \ref{stability}), we see that, for any $L>0$,
\[
\int_{-L}^{L}\left({1\over 2}|\dot \xi(s)|^2-V(\xi(s))-p\cdot \dot \xi(s)+\overline H(p)\right)\,ds=v(\xi(L))-v(\xi(-L)).
\]
Hence $\xi$ is a universal global characteristic associated with $p$.
\end{proof}

For smooth $V$, two different orbits in the same Aubry set cannot intersect and two different absolute minimizers of the same action cannot intersect twice.  However, both situations could happen with merely continuous $V$.  Consequently, the structure of orbits on $\mathcal {U}_p$ might be very complicated.    Below we provide two procedures to join different pieces of two global characteristics, which will be used later to select nice minimizing orbits and then simplify the topology of interacting curves.

\begin{rem}\label{glue} 
Given $p\in  \R^n$ and $v$ is a viscosity solution of \eqref{E-p},  suppose that $\xi_1$ and $\xi_2$ are two global characteristics associated with $v$ satisfying that  for some $t_1,t_2\in \R$, 
\[
\xi_1(t_1)=\xi_2(t_2).
\]
Define
\[
\xi_3(t)=
\begin{cases}
\xi_1(t)  \quad &\text{ for $t\leq t_1$}\\
\xi_2(t-t_1+t_2) \quad &\text{ for $t\geq t_1$}.
\end{cases}
\]
See Figure \ref{fig1}.
\begin{center}
\includegraphics[scale=0.5]{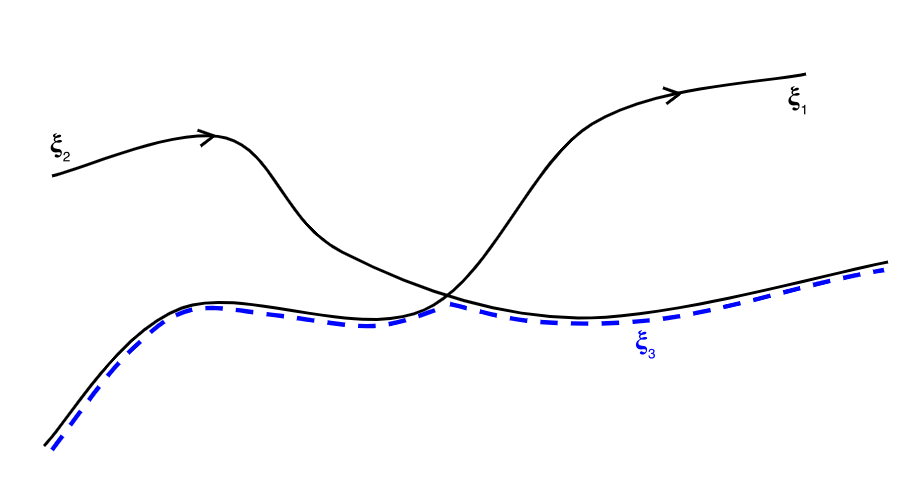}
\captionof{figure}{Formation of the curve $\xi_3$}\label{fig1}
\end{center}
Then, it is easy to see that $\xi_3$ is also a global characteristic associated with $v$. In fact,  for $a \leq t_1 \leq b$,
\[
v(b)-v(t_1)=\int_{t_1}^{b}\left({1\over 2}|\dot \xi_3(s)|^2-V(\xi_3(s))-p\cdot \dot \xi_3(s)+\overline H(p)\right)\,ds
\]
and
\[
v(t_1)-v(a)=\int_{a}^{t_1}\left({1\over 2}|\dot \xi_3(s)|^2-V(\xi_3(s))-p\cdot \dot \xi_3(s)+\overline H(p)\right)\,ds.
\]
Thus
\begin{align*}
v(b)-v(a)&=v(b)-v(t_1)+v(t_1)-v(a)\\
&=\int_{a}^{b}\left({1\over 2}|\dot \xi_3(s)|^2-V(\xi_3(s))-p\cdot \dot \xi_3(s)+\overline H(p)\right)\,ds.
\end{align*}

\end{rem}

\begin{rem}[Crossing of two universal global characteristics] \label{rem:two-curves}
Suppose that  $p,p'\in F_c$, and $\xi$ and $\tilde \xi$ are orbits in $\mathcal{U}_p$ and $\mathcal{U}_{p'}$, respectively. 
 Assume that there exist $t_1, t_2, t_{1}', t_{2}' \in \R$ such that,  for $i=1,2$, 
\[
P_i=\xi (t_i)=\tilde \xi(t_{i}').
\]
Now we present how to construct new orbits on  $\mathcal{U}_p$ and $\mathcal{U}_{p'}$ by joining different pieces of $\xi$ and $\tilde \xi$. Without loss of generality,  we assume that $t_1<t_2$. 
There are two cases.

\medskip

\noindent {\bf Case 1.} $t_{1}'<t_{2}'$.  
Define 
\[
\xi_2(t)=
\begin{cases}
\xi(t)  \quad &\text{ for $t\leq t_1$},\\
\tilde \xi(t+t_{1}'-t_1) \quad &\text{ for $t_1\leq t\leq t_1+t_{2}'-t_{1}'$},\\
\xi(t+t_2-t_1-t_{2}'+t_{1}') \quad &\text{ for $t_1+t_{2}'-t_{1}'\leq t$}.
\end{cases}
\]
See Figure \ref{fig2}.

\noindent {\bf Case 2.}  $t_{1}'>t_{2}'$.  Define 
\[
\xi_3(t)=
\begin{cases}
\xi(t)  \quad &\text{ for $t\leq t_1$},\\
\tilde \xi(t_1+t_{1}'-t) \quad &\text{ for $t_1\leq t\leq t_1+t_{1}'-t_{2}'$},\\
\xi(t+t_2-t_1-t_{1}'+t_{2}') \quad &\text{ for $t_1+t_{1}'-t_{2}'\leq t$}.
\end{cases}
\]

\begin{center}
\includegraphics[scale=0.5]{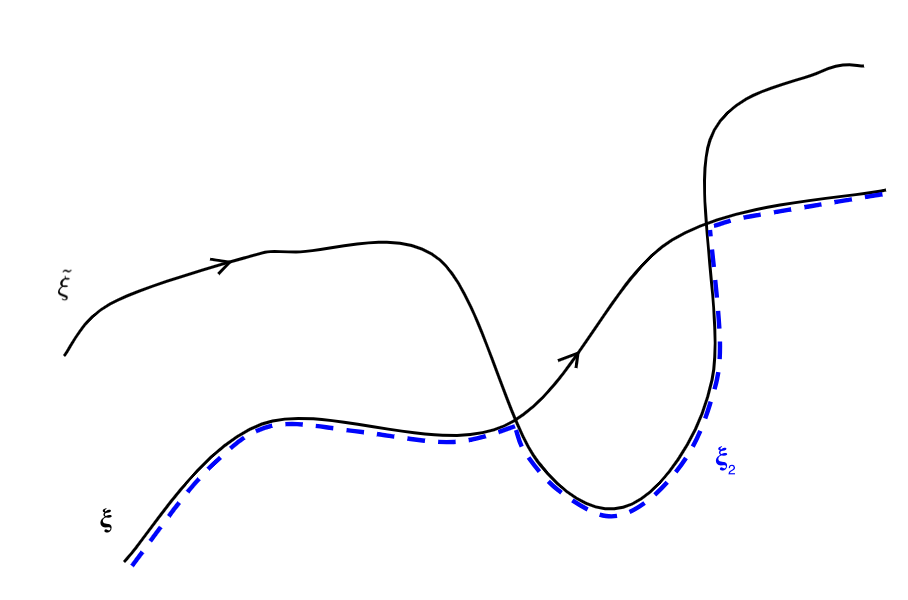}
\captionof{figure}{Combining two universal global characteristics}\label{fig2}
\end{center}

\end{rem}

Since both $\xi$ and $\tilde \xi$ are absolute action minimizing curves connecting $P_1$ and $P_2$, we have that 
\[
\int_{t_1}^{t_2}\left({1\over 2}|\dot \xi(s)|^2-V(\xi(s))+c\right)\,ds=\left |\int_{t_{1}'}^{t_{2}'}\left({1\over 2}|\dot {\tilde \xi}(s)|^2-V(\tilde \xi(s))+c\right)\,ds\right|. 
\]
Accordingly, we have  the following conclusion. 

\begin{cor}\label{crossing-join}
 The above $\xi_2$ or $\xi_3$ belongs to $\mathcal{U}_{p}$.
\end{cor}

\begin{defn}\label{adjustment} 
We say that $\xi_2$ or $\xi_3$ constructed in Remark \ref{rem:two-curves}  the adjustment of $\xi$ with respect to $\tilde \xi$ between $t_1$ and $t_2$. 
\end{defn}

\begin{lem}\label{direction} 
Let $p\in \R^n$, and $v$ be a viscosity solution of \eqref{E-p}. 
Let  $\xi:\R\to \R^n$ be a global characteristic associated with $v$. 
Assume that there exists $\{t_m\} \to \pm \infty$ such that ${\xi(t_m)/t_m}$ converges as $m\to\infty$.
Then, 
\[
\lim_{m\to \infty}{\xi(t_m)\over t_m}\in\partial \overline H(p).
\]
\end{lem}

\begin{proof}
Suppose that  
\[
\lim_{m\to \infty}{\xi(t_m)\over t_m}=q
\]
for some $q\in \R^n$.  
It suffices to show that 
\begin{equation}\label{support}
\ol H(p')\geq \ol H(p)+q\cdot (p'-p) \quad \text{ for all $p'\in \R^n$.}
\end{equation}
Let $v'$ be a viscosity solution of \eqref{E-p} with $p=p'$.  
By Lemma \ref{calib},
\[
p'\cdot \xi(t_m)+v'(\xi(t_m))-(p'\cdot \xi(0)+v'(\xi(0)))\leq \int_{0}^{t_m}\left({1\over 2}|\dot \xi(s)|^2-V(\xi(s))+\overline H(p')\right)\,ds.
\]
Meanwhile,  
\[
p\cdot \xi(t_m)+v(\xi(t_m))-(p\cdot \xi(0)+v(\xi(0)) = \int_{0}^{t_m}\left({1\over 2}|\dot \xi(s)|^2-V(\xi(s))+\overline H(p)\right)\,ds.
\]
Taking the difference of the above two equations, dividing both sides by $t_m$, and sending $m\to \infty$, we derive \eqref{support}.
\end{proof}

\begin{defn}  
For $p\in  \R^n$ and $v$ is a viscosity solution of \eqref{E-p},  a global characteristic  $\xi:\R\to \R^n$ associated with $v$ is called  periodic if there exist $T>0$ and $q\in \Z^n$  such that 
\[
\xi(t+T)-\xi(t)=q  \quad \text{ for all $t\in \R$.}
\]
In this case, $q/T$ is called the rotation vector of $\xi$. 
\end{defn}

Owing to Lemma \ref{direction},  the rotation vector 
\[
{q\over T}\in \partial \overline H(p).
\]
Also, it is clear that every periodic  global characteristic associated with some $v$ must be a universal global characteristic.

\begin{cor}\label{two-integer-points} 
For $p\in  \R^n$ and $v$ is a viscosity solution of \eqref{E-p},  let $\xi:\R \to \R^n$ be a global characteristic associated with $v$.  
If there exist $t_1<t_2$ such that 
\[
\xi(t_2)-\xi(t_1)=q\in \Z^n,
\]
 then
\[
{q\over {t_2-t_1}}\in  \partial \overline H(p).
\]
\end{cor}

\begin{proof}
Let $T=t_2-t_1$.  
Define $\tilde \xi:\R\to \R^n$ as
\[
\tilde \xi(t)=\xi (t-kT+t_1)+kq  \quad \text{if $t\in  [kT, (k+1)T]$},
\]
for all $k\in \Z$.  Then $\tilde \xi(0)=\xi(t_1)$ and $\tilde \xi(t+T)=\tilde \xi(t)+q$ for all $t\in \R$.  Since 
\[
0=v(\xi(t_2))-v(\xi(t_1))=\int_{t_1}^{t_2}\left({1\over 2}|\dot \xi(t)|^2-V(\xi(t))-p\cdot \dot \xi(t) + \ol H(p) \right)\,dt,
\]
for any $m\in \N$, we have that
\begin{align*}
I(\tilde \xi, [-mT,mT])&=\sum_{k=-m}^{m-1}I(\tilde \xi, [kT, (k+1)T])\\
&=\sum_{k=-m}^{m-1}I(\tilde \xi, [0, T])=\sum_{k=-m}^{m-1}I(\xi, [t_1, t_2])\\
&=0=v(\tilde \xi(mT))-v(\tilde \xi(-mT)).
\end{align*}
Here,
\[
I(\eta, [a,b])=\int_{a}^{b}\left({1\over 2}|\dot \eta(t)|^2-V(\eta(t))-p\cdot \dot \eta(t) + \ol H(p) \right)\,dt.
\]
Then, $\tilde \xi$ is a periodic global characteristic associated with $v$. 
In fact, $\tilde \xi\in \mathcal{U}_p$.  
The conclusion follows from Lemma \ref{direction}. 
\end{proof}

\begin{lem}\label{p-orbit}
When $n=2$,  for every  $q\in \Z^2$ and $c>\max_{\Bbb T^n}V$, there exists $p_q\in F_c$ such that $\mathcal{U}_{p_q}$ has a periodic orbit $\xi$ such that, for some $T>0$,
\[
\xi(t+T)-\xi(t)=q \quad \text{ for all $t\in \R$}.
\]
\end{lem}

\begin{proof}
This result is well-known for smooth $V$.  
Similar to the proof of the stability of global characteristics in Lemma \ref{stability},  the merely continuous version can be established by approximating $V$ with smooth  periodic functions under the maximum norm.
\end{proof}

\begin{lem}\label{equal-aubry} 
Suppose that there exist $p_0,p_1\in \R^n$ and  $\lambda\in (0,1)$ such that, for $p_\lambda=\lambda p_0+(1-\lambda)p_1$,
\[
\overline H(p_\lambda)=\lambda \overline H(p_0)+ (1-\lambda) \overline H(p_1).
\]
Then,
\[
\overline H(p_\lambda)=\overline H(p_0)=\overline H(p_1),\qquad {\mathcal{A}}_{p_\lambda}\subset  {\mathcal{A}}_{p_0}\cap{\mathcal{A}}_{p_1},
\]
and
\[
{\mathcal{U}}_{p_\lambda}\cap {\mathcal{U}}_{p_0}\cap{\mathcal{U}}_{p_1}\not= \emptyset. 
\]
\end{lem}

\begin{proof}
We divide the proof into two steps.

\medskip

\noindent{\bf Step 1.}   
For $x\in {\mathcal{A}}_{p_\lambda}$,   there exist $\{t_m\}\to \infty$ and a sequence of curves $\gamma_m: [0,t_m]\to \R^n$ such that $\gamma_m(0)=x$, $\gamma_m(t_m)\in x+\Z^n$, and
\[
\lim_{m\to \infty}\int_{0}^{t_m}\left({1\over 2}|\dot \gamma_m(s)|^2-V(\gamma_m(s))-p_{\lambda}\cdot \dot \gamma_m(s)+\ol H(p_\lambda)\right)\,ds=0.
\]
Let
\[
A_m=\lambda\int_{0}^{t_m}\left({1\over 2}|\dot \gamma_m(s)|^2-V(\gamma_m(s))-p_0\cdot \dot \gamma_m(s)+\ol H(p_0)\right)\,ds,
\]
and
\[
B_m=(1-\lambda)\int_{0}^{t_m}\left({1\over 2}|\dot \gamma_m(s)|^2-V(\gamma_m(s))-p_1\cdot \dot \gamma_m(s)+\ol H(p_1)\right)\,ds.
\]
Then
\[
\int_{0}^{t_m}\left({1\over 2}|\dot \gamma_m(s)|^2-V(\gamma_m(s))-p_{\lambda}\cdot \dot \gamma_m(s)+\ol H(p_\lambda)\right)\,ds=A_m+B_m.
\]
By Lemma \ref{calib},  $A_m, B_m\geq 0$. 
Hence, 
\[
\lim_{m\to \infty}A_m=\lim_{m\to \infty}B_m=0.
\]
Meanwhile, by the definition of $G_p(x,x)$, it is obvious that 
\[
0\leq \lambda G_{p_0}(x,x)\leq \lim_{m\to \infty}A_m \quad \text{ and } \quad 0\leq (1-\lambda)G_{p_1}(x,x)\leq \lim_{m\to \infty}B_m.
\]
Accordingly, $G_{p_0}(x,x)=G_{p_1}(x,x)=0$. Then $x\in \mathcal{A}_{p_0}\cap \mathcal{A}_{p_1}$.  This implies that 
\[
 {\mathcal{A}}_{p_\lambda}\subset {\mathcal{A}}_{p_0}\cap {\mathcal{A}}_{p_1}.
\]

In addition,  as in the proof of Lemma \ref{Aubry-orbit},  a suitable reparametrization of $\{\gamma_m\}$ gives a sequence of curves that converges to a common orbit in ${\mathcal{U}}_{p_\lambda}\cap {\mathcal{U}}_{p_0}\cap{\mathcal{U}}_{p_1}$.

\medskip

\noindent{\bf Step 2.} 
Choose an orbit $\xi \in {\mathcal{U}}_{p_0}\cap{\mathcal{U}}_{p_1}$. 
Then, by (\ref{constant-energy-aubry}), 
\[
{1\over 2}|\dot \xi(t)|^2+V(\xi(t))=\ol H(p_0)  \quad \text{for a.e $t\in \R$,}
\]
and
\[
{1\over 2}|\dot \xi(t)|^2+V(\xi(t))=\ol H(p_1)  \quad \text{for a.e $t\in \R$}.
\]
Therefore, $\ol H(p_0)=\ol H(p_1)$. 

\end{proof}


\section{Proof of Theorem \ref{thm:main}}\label{sec:proof}
We are in the setting of Theorem \ref{thm:main} in this section.  
Part (1) follows immediately from Lemma \ref{equal-aubry}. 
We now prove part (2).  Throughout this section, we assume that $c>\max_{\Bbb T^2}V$.  
Hence any orbit on $\mathcal{U}_p$ for $p\in F_c$ does not intersect with itself since $u=p\cdot x+v$ is strictly increasing along any orbit. 

We argue by contradiction.  
Suppose that $F_c$ contains a line segment of an irrational slope.  
Assume that $p_0$ and $p_1$ are two points in the interior of the line segment. 
According to Lemma \ref{equal-aubry}, 
\[
\mathcal{A}_{p_0}= \mathcal{A}_{p_1},
\]
and
\[
\text{$p_0-p_1$ is an irrational vector}. 
\]
Then, the  outward unit normal vector $\vec{n}$ is also irrational, and
\[
\partial \overline H(p_0)=\partial \overline H(p_1)=\left\{\lambda \vec{n} \,:\,  \lambda\in [\alpha, \beta]\right\}
\]
for two positive numbers $0<\alpha<\beta$.  
Without loss of generality, we assume that 
\begin{equation}\label{h-direction}
\vec{n}\cdot (1,0)>0.
\end{equation}
Let  $v_0$ and $v_1$ be viscosity solutions to  \eqref{E-p} corresponding to $p=p_0$ and $p=p_1$, respectively.  
Write
\[
\mathcal{U}=\mathcal{U}_{p_0}\cap \mathcal{U}_{p_1},
\]
and
\[
S=\bigcup_{\xi \in \mathcal{U}}\xi(\R)\subset \R^2.
\]
Owing to Lemma \ref{equal-aubry}, $\mathcal{U}\not=\emptyset$.  Also,  by Lemma \ref{regularity}, $u_0-u_1$ is differentiable at $x\in S$, and 
\begin{equation}\label{zero-gradient}
D(u_0-u_1)(x)=0 \quad \text{ for $x\in S$}. 
\end{equation}
Here, $u_0(x)=p_0\cdot x+v_0(x)$, and $u_1(x)=p_1\cdot x+v_1(x)$ for $x\in \R^2$. 

\medskip

Thanks to Lemma \ref{p-orbit},  we can choose $p'\in F_c$ such that $\mathcal {U}_{p'}$ contains a periodic orbit $\eta$ such that for some $T>0$
\[
\eta(t+T)-\eta(t)=(0,1)=e_2 \quad \text{ for all $t\in \R$}.
\]
For $e_1=(1,0)$, denote
\begin{equation}\label{h-distance}
\Lambda=\max\{|e_1\cdot (x-y)|\,:\,  x,y\in \eta(\R)\}.
\end{equation}
Choose a positive integer $J>\Lambda+1$ and for $k\in \Z$, denote 
\[
\eta_k=\eta+k(J,0).
\]
\begin{center}
\includegraphics[scale=0.5]{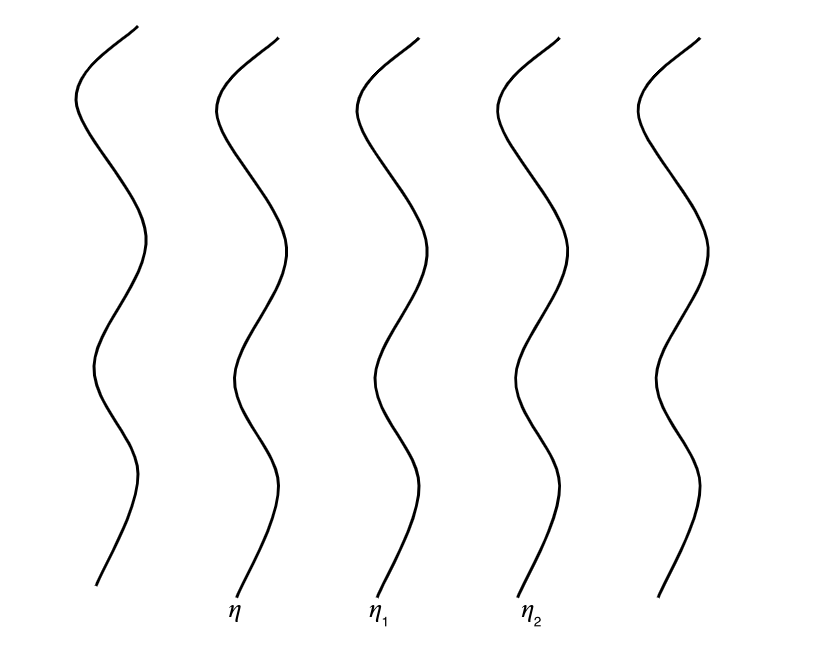}
\captionof{figure}{Family of  $\{\eta_k\}_{k\in \Z}$}\label{fig3}
\end{center}
Clearly,  these curves are mutually disjoint.  
See Figure \ref{fig3}.

Throughout this section, $\xi$ represents a given orbit on $\mathcal{U}$.  
Then, owing to Lemma \ref{direction},
\begin{equation}\label{normal-direction}
\lim_{t\to \infty}{\xi(t)\over |\xi(t)|}=\lim_{t\to -\infty}{-\xi(t)\over |\xi(t)|}=\vec{n}.
\end{equation}
Then, $\xi$ intersects each $\eta_k$. 

For   $k\in \Z$, write
\[
t_{k,+}=\max\{t\in \R\,:\,  \xi(t)\in  \eta_k(\R)\} \quad \text{ and } \quad t_{k,-}=\min\{t\in \R\,:\,  \xi(t)\in  \eta_k(\R)\}.
\]
Owing to (\ref{normal-direction}), both $t_{k,+}$ and $t_{k,-}$ are finite.  
See Figure \ref{fig4}.
\begin{center}
\includegraphics[scale=0.5]{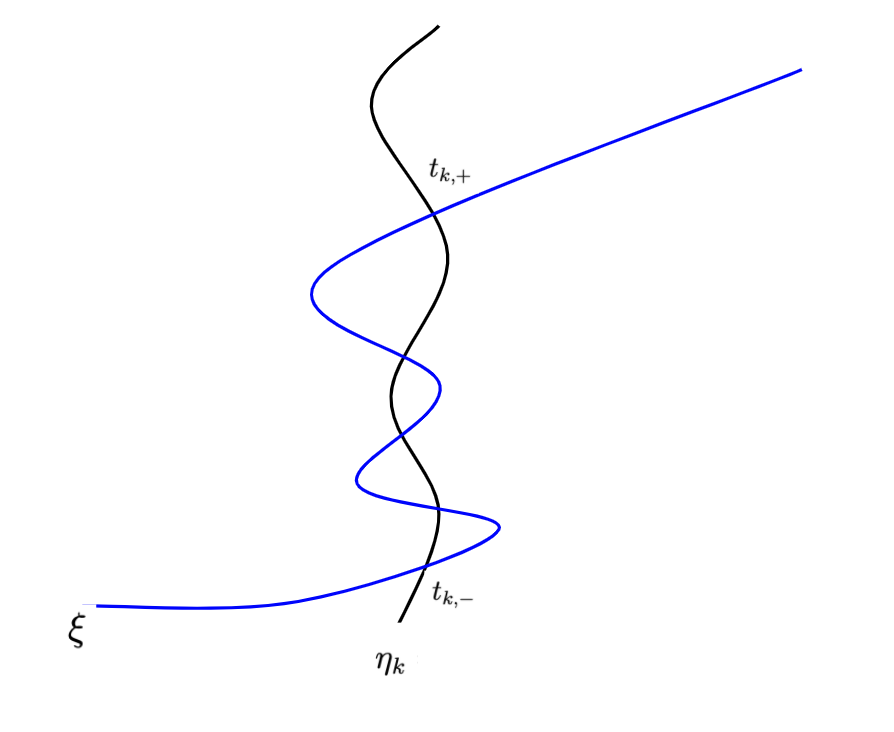}
\captionof{figure}{Intersections of $\eta_k$ and $\xi$}\label{fig4}
\end{center}

\begin{lem}\label{monotonicity} 
For $k\in \Z$, assume that 
\[
\xi(t_{k,+})=\eta_k(\theta_+)   \quad \text{ and } \quad \xi(t_{k,-})=\eta_k(\theta_-)
\]
for $\theta_-,\theta_+ \in \R$.
Then,
\begin{itemize}
\item[(1)]
\[
\xi(\R)\cap \eta_k(\R)\subset \left\{\eta_k(t)\,:\,  \min\{\theta_{+},  \theta_{-}\}\leq   t\leq \max\{\theta_{+},  \theta_{-}\}\right\};
\]

\item[(2)]
\[ |\theta_+-\theta_-|<T.
\]
\end{itemize}
\end{lem}

\begin{proof}
We first prove (1).   
Were the conclusion of (1)  not true, there would exist $t_0\in (t_{k,-}, t_{k,+})$ such that 
\[
\xi(t_0)=\eta_k(\theta)
\]
for some $\theta<\min\{\theta_{+},  \theta_{-}\}$ or $\theta> \max\{\theta_{+},  \theta_{-}\}$.  
Without loss of generality, we assume that $\theta>\theta_{+}>\theta_{-}$. 
Then,
\[
\underbrace{\int_{\theta_{-}}^{\theta_{+}}\left({1\over 2}|\dot \eta_k(s)|^2-V(\eta_k(s))+c\right)\,ds}_{A}=\underbrace{\int_{t_{k,-}}^{t_{k,+}}\left({1\over 2}|\dot \xi(s)|^2-V(\xi(s))+c\right)\,ds}_{B}
\]
as both $\eta_k$ and $\xi$ are absolute minimizers of the action connecting $\eta_k(\theta_-)$ and $\eta_k(\theta_+)$. 
Also, 
\[
\underbrace{\int_{\theta_{-}}^{\theta}\left({1\over 2}|\dot \eta_k(s)|^2-V(\eta_k(s))+c\right)\,ds}_{C}=\underbrace{\int_{t_{k,-}}^{t_{0}}\left({1\over 2}|\dot \xi(s)|^2-V(\xi(s))+c\right)\,ds}_{D}
\]
since both $\eta_k$ and $\xi$ are absolute minimizers of the action connecting $\eta_k(\theta_-)$ and $\eta_k(\theta)$. On the other hand, it is obvious that 
\[
B> D  \quad \text{ and } \quad A<C.
\]
This is a contradiction.

\smallskip

Next we prove (2).  
Again we argue by contradiction.  
Without loss of generality,  assume that  $\theta_+-\theta_-\geq T$. 
See Figure \ref{fig5}.
\begin{center}
\includegraphics[scale=0.5]{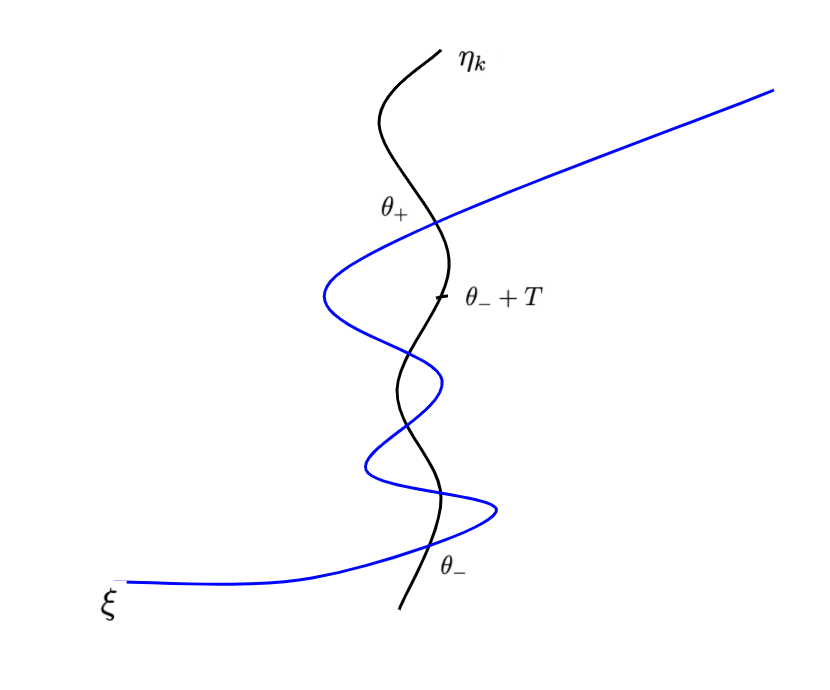}
\captionof{figure}{The situation where $\theta_+-\theta_-\geq T$}\label{fig5}
\end{center}

Let $\xi_1$ be the adjustment of $\xi$ with respect to $\eta_k$ between $t_{k,-}$ and $t_{k,+}$ (see Definition \ref{adjustment}).  
Then,  owing to (1),  $\xi_1(\R)\cap \eta_k(\R)=\eta_k([\theta_-, \theta_+])$. 
In particular, 
 $\xi_1$ would contain two points $A=\eta_k(\theta_-)$ and $B=\eta_k(\theta_-+T)=A+(0,1)$. 
Owing to Corollary  \ref{two-integer-points}, $\vec{n}$, the normal vector of $F_c$ at $p_0$, is parallel to $(0,1)$, which contradicts  the assumption that $\vec{n}$ is irrational.
\end{proof}

\begin{rem} 
(1) in the the above lemma actually says that the intersection parameters are monotonic.  Precisely speaking, for $t_1<t_2<t_3$, if 
\[
\xi(t_i)=\eta_k(\theta_i) \quad \text{  for $i=1,2,3$,}
\]
then  we have either  $\theta_1<\theta_2<\theta_3$ or $\theta_1>\theta_2>\theta_3$.
\end{rem}

Write
\[
L_A=\int_{0}^{T}\left({1\over 2}|\dot \eta(s)|^2-V(\eta(s))+c\right)\,ds,
\]
which is exactly the action of one cycle of $\eta$. 

Owing to \eqref{h-direction} and \eqref{normal-direction},  $\xi$ will intersect $\eta_k$ before it intersects  $\eta_{k+1}$.   The following lemma says that $\xi$ will not intersect $\eta_k$ again after it intersects $\eta_{k+1}$ for large enough $J$.

\begin{lem}\label{intersection-order} 
Assume that 
\[
J>\max\left\{\Lambda+1, \ {L_A\over \sqrt{c-\max_{\T^2}V}}+\Lambda\right \}.
\]
Then, 
\[
t_{k,+}<t_{k+1,-} \quad \text{ for all $k\in  \Z$}.
\]
Recall that $\Lambda=\max\{|e_1\cdot (x-y)|\,:\,  x,y\in \eta(\R)\}$. 
\end{lem}

\begin{proof}
If the conclusion of the lemma were false, then we would have  
\[
t_{k+1,-}\in (t_{k,-}, t_{k,+}).
\] 
See Figure \ref{fig6}.
\begin{center}
\includegraphics[scale=0.5]{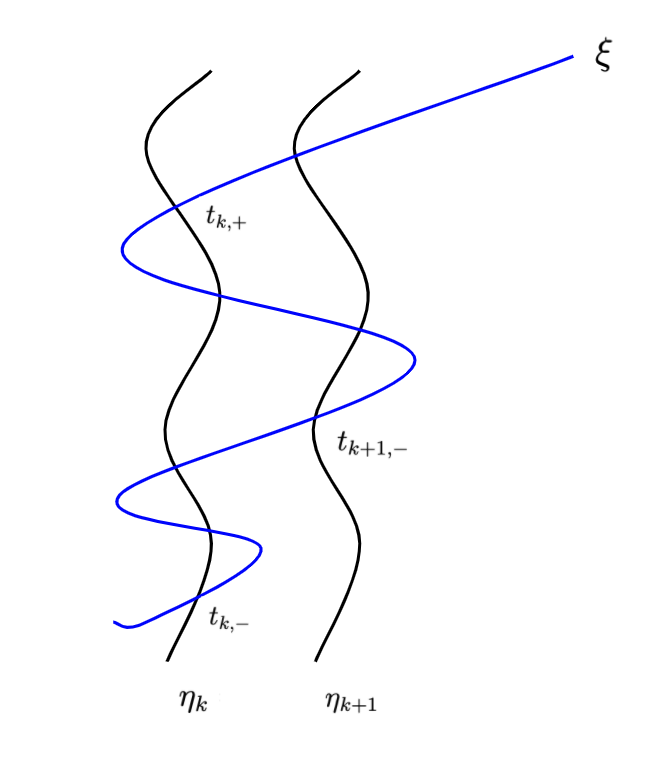}
\captionof{figure}{The situation where $t_{k+1,-}\in (t_{k,-}, t_{k,+})$}\label{fig6}
\end{center}
Note first that
\[
2(J-\Lambda)\leq |\xi(t_{k+1,-})-\xi(t_{k,-})|+|\xi(t_{k+1,-})-\xi(t_{k,+})|\leq \int_{t_{k,-}}^{ t_{k,+}}|\dot \xi(s)|\,ds.
\]
Meanwhile, 
\begin{align*}
 \int_{t_{k,-}}^{ t_{k,+}}\left({1\over 2}|\dot \xi(s)|^2-V(\xi(s))+c\right)\,ds&\geq \sqrt{2}\int_{t_{k,-}}^{ t_{k,+}}\sqrt{c-V(\xi(s))}|\dot \xi(s)|\,ds\\
&\geq \sqrt{c-\max_{\Bbb T^2}V}\int_{t_{k,-}}^{ t_{k,+}}|\dot \xi(s)|\,ds.
\end{align*}
Also, 
\[
\int_{t_{k,-}}^{ t_{k,+}}\left({1\over 2}|\dot \xi(s)|^2-V(\xi(s))+c\right)\,ds=\left|\int_{\theta_{k,-}}^{ \theta_{k,+}}\left({1\over 2}|\dot \eta_k(s)|^2-V(\eta_k(s))+c\right)\,ds\right|\leq L_A.
\]
Here, $\xi(t_{k,+})=\eta_k(\theta_{k,+})$,  and $\quad \xi(t_{k,-})=\eta_k(\theta_{k,-})$. 
The inequality in the above is due to (2) in  Lemma \ref{monotonicity}.  The absolute sign $|\cdot |$ is added because  $\theta_{k,-}$ might be bigger than $\theta_{k,+}$.
Therefore,
\[
2(J-\Lambda)\leq  {L_A\over \sqrt{c-\max_{\T^2}V}},
\]
which contradicts the choice of $J$.
\end{proof}

For two orbits $\xi_1$ and $\xi_2$ on $\mathcal{U}$, we define the distance
\[
d(\xi_1,\xi_2)=\sum_{k\in \Z}{\arctan (d(\xi_1, \xi_2,k))\over |k|^2+1}.
\]
Here,
\[
d(\xi_1, \xi_2,k)=\min\{|\theta_1-\theta_2| \,:\,  \theta_1\in A_{1,k},\ \theta_2\in A_{2,k}\},
\]
where, for $i=1,2$, 
\[
A_{i,k}=\{\theta\in \R \,:\,  \eta_k(\theta)\in  \xi_i(\R)\}.
\]
Clearly,   the distance function $d(\cdot, \cdot, k)$ (hence $d(\cdot, \cdot)$) is lower semicontinuous with respect to orbits on $\mathcal{U}$, i.e.,  if $\xi_{i,n}\to \xi_i$ locally uniformly for $i=1,2$, then
\[
\liminf_{m\to \infty}d(\xi_{1,n}, \xi_{2,n})\geq d(\xi_1, \xi_2).
\]

Define
\begin{equation}\label{def-intersection}
\mathcal{I}=\{\theta \in \R \,:\,  \text{there exists an orbit $\xi\in \mathcal{U}$ such that $\eta(\theta)\in \xi(\R)$}\}.
\end{equation}
Due to the $T$-periodicity of $\eta$ and the fact that  the set $\mathcal{U}$ is closed under limits of orbits,   $\mathcal{I}$ is a $T$-periodic closed set.  

To finish the proof,  the next step  is to show that  for $u_0(x)=p_0\cdot x+v_0(x)$, and $u_1(x)=p_1\cdot x+v_1(x)$ for $x\in \R^2$, 
\begin{equation}\label{eq:main-con}
\text{$u_0-u_1$ is constant on $\mathcal{I}$},
\end{equation}
which will lead to a contradiction since $\lim_{\theta\to \infty} |u_0(\eta(\theta))-u_1(\eta(\theta))|=\infty$. 
Suppose that 
\[
\R\backslash \mathcal{I}=\bigcup_{i=1}^{\infty}\,(a_i,b_i),
\]
where $\{(a_i,b_i)\}_{i\geq 1}$ are disjoint open intervals. Obviously, $(a_i, b_i)\subseteq (a_i, a_i+T)$ for each $i\in \N$ since $\xi\in \mathcal{U}\Rightarrow \xi+(0,1)\in \mathcal{U}$. 

\begin{lem} \label{lem:ee}
For all $j\in \N$, 
\begin{equation}\label{equal-end}
u_0(\eta(a_j))-u_1(\eta(a_j))=u_0(\eta(b_j))-u_1(\eta(b_j)).
\end{equation}

\end{lem}

\begin{proof}
We only need to prove the claim for $j=1$.
Owing to the lower semicontinuity of the distance function, we may choose $\xi_1$ and $\xi_2$ in $\mathcal{U}$ such that

\begin{itemize}
\item[(i)] $\xi_1(0)=\eta(a_1)$ and $\xi_2(0)=\eta(b_1)$.

\item[(ii)] $d(\xi_1,\xi_2)$ attains the minimum value among all curves in $\mathcal{U}$ satisfying (i).
\end{itemize}
To simplify the associated topology between curves,  we adjust $\xi_1$ and $\xi_2$ with respect to each $\eta_k$ between $t_{i,k,-}$ and $t_{i, k,+}$ for $i=1,2$ respectively (see Definition \ref{adjustment}). Here 
\[
t_{i, k,+}=\max\{t\in \R\,:\,  \xi_i(t)\in  \eta_k(\R)\} \quad \text{ and } \quad t_{i, k,-}=\min\{t\in \R\,:\,  \xi_i(t)\in  \eta_k(\R)\}.
\]
By Lemma \ref{intersection-order}, for $i=1,2$,  the time intervals $(t_{i, k,-}, t_{i, k,+})$ are mutually disjoint, i.e., 
\[
  ... <t_{i,-1,-}\leq t_{i,-1,+}<t_{i,0,-}\leq t_{i,0,+}<t_{i,1,-}\leq t_{i,1,+}<t_{i,2,-}\leq t_{i,2,+}<..
\]
Hence the adjustments are well defined.   

In addition, thanks to (1) of Lemma \ref{monotonicity}, the distance between two adjusted orbits is not greater than  $d(\xi_1,\xi_2)$.   Thus two adjusted orbits also satisfy (i)-(ii) above.  
By abuse of notations,  we still use $\xi_1$ and $\xi_2$ to represent corresponding adjusted orbits. 
Accordingly, we may assume that for each $k\in \Z$, and $i=1,2$,
\begin{equation}\label{good-topology}
\xi_i(\R)\cap \eta_k(\R)=\xi_i([t_{i,k,-},\ t_{i,k,+}])=\eta_k([\theta_{i,k,-},\ \theta_{i,k,+}]). 
\end{equation}
Here $\xi_i(t_{i,k,+})=\eta_k(\theta_{i, k,+})$  and $\quad \xi(t_{i, k,-})=\eta_k(\theta_{i, k,-})$.  
It could happen that $\theta_{i,k,+}<\theta_{i, k,-}$.  
In terms of topology, the above basically  plays the role like that  $\xi_i$ and $\eta_k$ only intersect once for smooth $V$.
We consider two cases.

\medskip

\noindent {\bf Case 1.}  $\xi_1(\R)\cap \xi_2(\R)\not=\emptyset$. Assume that $\xi_1(t_1)=\xi_2(t_2)$ for some $t_1, t_2\in \R$. 
See Figure \ref{fig7}.
\begin{center}
\includegraphics[scale=0.5]{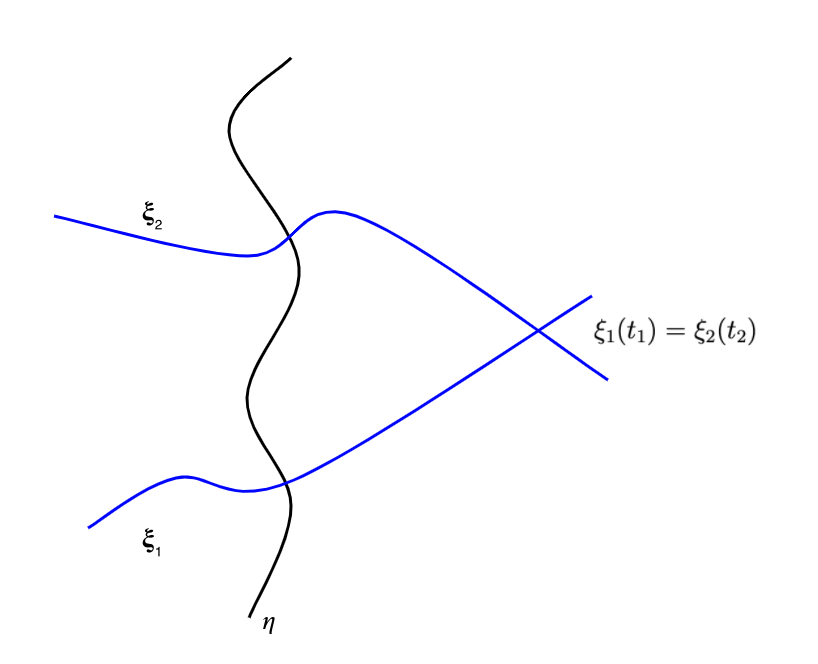}
\captionof{figure}{The situation where $\xi_1(\R)\cap \xi_2(\R)\not=\emptyset$}\label{fig7}
\end{center}
Then,
\[
u_0(\xi_1(t_1))-u_0(\xi_1(0))=u_1(\xi_1(t_1))-u_1(\xi_1(0))=\int_{0}^{t_1}\left({1\over 2}|\dot \xi_1(s)|^2-V(\xi_1(s))+c\right)\,ds,
\]
and
\[
u_0(\xi_2(t_2))-u_0(\xi_2(0))=u_1(\xi_2(t_2))-u_1(\xi_2(0))=\int_{0}^{t_2}\left({1\over 2}|\dot \xi_2(s)|^2-V(\xi_2(s))+c\right)\,ds.
\]
Taking the difference of the two equations above leads to the claim. 

\medskip

\noindent {\bf Case 2.} $\xi_1(\R)\cap \xi_2(\R)=\emptyset$.  For each $k\in \Z$, let 
\[
d_1(k)=\max\{\theta \,:\,  \eta_k(\theta)\in  \xi_1(\R)\cap \eta_k(\R)\},
\]
and
\[
d_2(k)=\min\{\theta \,:\,  \eta_k(\theta)\in  \xi_2(\R)\cap \eta_k(\R)\}.
\]
Owing to two dimensional topology and (\ref{good-topology}),  we have that 
\[
d_1(0)=a_1  \quad \text{ and } \quad d_2(0)=b_1,
\]
and
\[
d_1(k)<d_2(k)  \quad \text{ for all $k\in \Z$}. 
\]
See Figure \ref{fig8}.
\begin{center}
\includegraphics[scale=0.5]{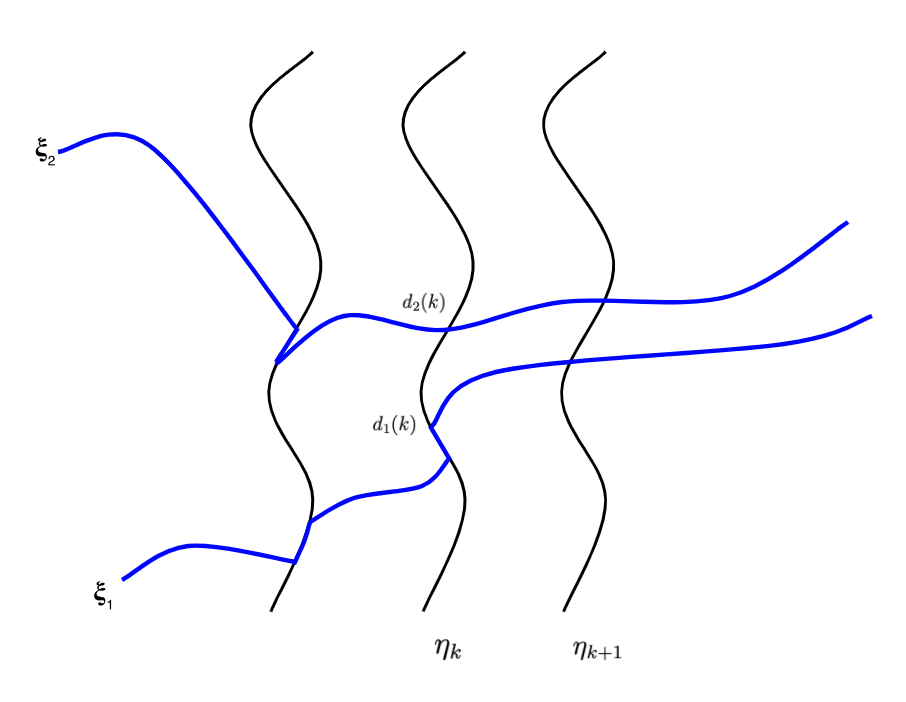}
\captionof{figure}{Positions of $d_1(k),d_2(k)$ }\label{fig8}
\end{center}
We claim that
\begin{equation}\label{eq:sum-k}
\sum_{k\in \Z}(d_2(k)-d_1(k))\leq T,
\end{equation}
which will be proved in Lemma \ref{lem:sum-k} below.
In particular, \eqref{eq:sum-k} implies 
\[
\lim_{k\to \infty} (d_2(k)-d_1(k))=0.
\]
Similar to Case 1 above, since $\eta_k(d_i(k))\in \xi_i(\R)$ for $i=1,2$, 
\[
u_0(\eta_k(d_1(k)))-u_0(\xi_1(0))=u_1(\eta_k(d_1(k)))-u_1(\xi_1(0)),
\]
and
\[
u_0(\eta_k(d_2(k)))-u_0(\xi_2(0))=u_1(\eta_k(d_2(k)))-u_1(\xi_2(0)). 
\]
Taking the difference of the two equations and sending $k\to \infty$, we obtain the claim.   

\end{proof}

\begin{proof}[Proof of Theorem \ref{thm:main}]
To finish the proof of the main result, we only need to prove \eqref{eq:main-con}.
We proceed by using \eqref{equal-end}.

Write
\[
g(t)=u_{0}(\eta(t))-u_{1}(\eta(t)).
\]
Then, $g(t)$ is Lipschitz continuous and by \eqref{zero-gradient}, 
\[
g'(t)=0  \quad \text{ for $t\in \mathcal{I}$}.
\]
See \eqref{def-intersection} for the definition of $\mathcal{I}$. Moreover, due to \eqref{equal-end},  for each $i\in \N$,
\[
g(a_i)=g(b_i).
\]
Define a new function
\[
h(t)=
\begin{cases}
g(t) \quad &\text{ for $t\in \mathcal{I}$}\\
g(a_i) \quad  &\text{ for $t\in (a_i, b_i)$, and $i\in \N$}.
\end{cases}
\]
Clearly, $h(t)$ is Lipschitz continuous and $h'(t)=0$ for a.e. $t\in \R$. 
Hence $h \equiv c$ for some constant $c\in \R$.  
This leads to
\[
g(t)= c  \quad \text{ for $t\in \mathcal{I}$}.
\]
This is absurd as for $m\in \N$, 
\[
|g(mT)-g(0)|\geq m|(p_0-p_1)\cdot (0,1)|-\max_{\T^2}|v_0|-\max_{\T^2}|v_1|,
\] 
which leads to $\lim_{t\to \infty}|g(t)|=+\infty$. 
\end{proof}

Finally, we give a proof of  \eqref{eq:sum-k}.
\begin{lem}\label{lem:sum-k}
In  Case 2 in the proof of Lemma \ref{lem:ee}, 
\begin{equation*}
\sum_{k\in \Z}(d_2(k)-d_1(k))\leq T.
\end{equation*}
\end{lem}

\begin{proof}
We first show that,  for all $k\in \Z$,  the open interval
\begin{equation}\label{claim-d-1-2}
(d_1(k),\  d_2(k))\subset \R\backslash \mathcal{I},
\end{equation}
that is, it is one of those open intervals $\{(a_j,b_j)\}_{j\geq 1}$. 
It suffices to show this for $k>0$ as the proof for $k<0$ is similar.   
We argue by contradiction.  
If this were not true, then there would exist $k\in \N$ and  $\tilde \xi\in \mathcal{U}$ such that 
\begin{equation}\label{0-position}
\tilde \xi(0)\in \{\eta_k(\theta)\,:\, \theta\in (d_1(k),\  d_2(k))\}.
\end{equation}
Since $\tilde \xi$ cannot pass the portion of $\eta$ on  $(a_0,b_0)$,  we deduce that, if trace backward along $\tilde \xi$,   it must intersect $\xi_1$ or $\xi_2$ before it  intersects $\eta$. 
See Figure \ref{fig9}.
\begin{center}
\includegraphics[scale=0.5]{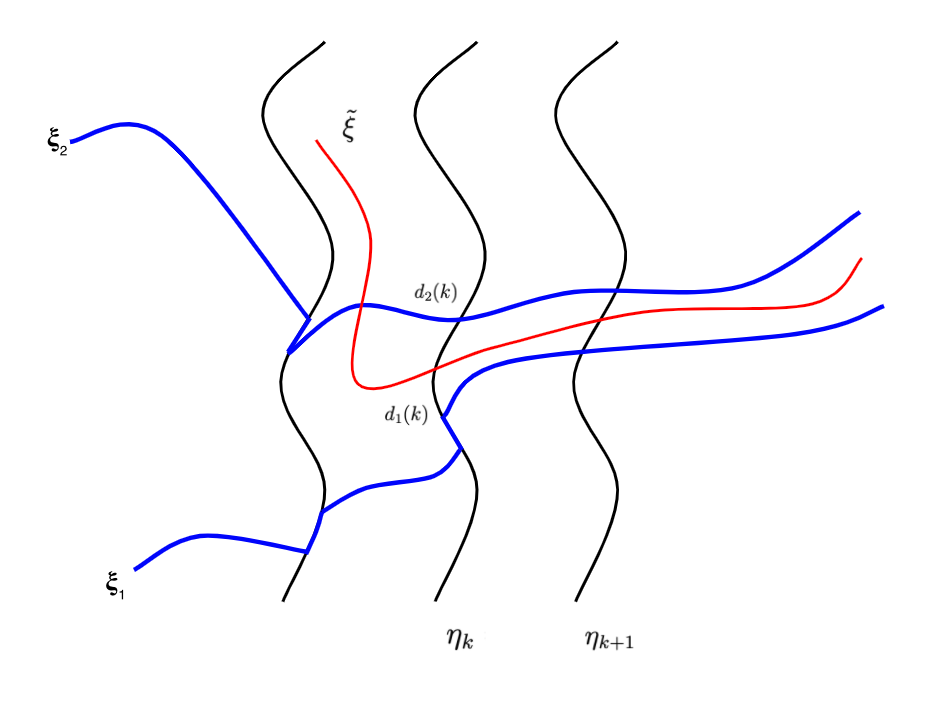}
\captionof{figure}{Relative position of $\tilde \xi$}\label{fig9}
\end{center}

Let
\[
t_{-}=\max\{t\leq 0 \,:\,  \tilde \xi(t)\in  \xi_1(\R)\cup\xi_2(\R)\}, \quad t_{+}=\inf\{t\geq 0\,:\,  \tilde \xi(t)\in  \xi_1(\R)\cup\xi_2(\R)\}.
\]
Then $t_{-}<0$ and  $t_+>0$. Note that $t_{+}$ could be $+\infty$. 

Now we will use the gluing property of Remark \ref{glue} to construct a new orbit in $\mathcal{U}$.  By two dimensional topology and \eqref{good-topology}, it is easy to see that, for each $k\in \Z$,
\begin{equation}\label{right-interval}
\tilde \xi((t_{-}, t_{+}))\cap (\eta_k(\R))\subset  \{\eta_k(\theta)\,:\, \theta\in (d_1(k),\  d_2(k))\}.
\end{equation}

Without loss of generality, we assume that $\tilde \xi(t_{-})\in \xi_2(\R)$ and $\tilde \xi(t_{+})\in \xi_j(\R)$ for $j=1$ or $j=2$ if $t_+<+\infty$.  Suppose that 
\[
\tilde \xi(t_{-})=\xi_2(\bar t_{-})  \quad \mathrm{and} \quad \tilde \xi(t_{+})=\xi_j(\bar t_{+}). 
\]
for $0<\bar t_{-}<\bar t_{+}$. Here $\bar t_{+}=+\infty$ if $ t_{+}=+\infty$. 
Let
\[
\xi_3(t)=
\begin{cases}
\xi_2(t) \quad &\text{ for $t\leq \bar t_{-}$},\\
\tilde \xi(t+t_{-}-\bar t_{-}) \quad &\text{ for $ \bar t_{-}\leq t\leq  t_{+}+\bar t_{-}-t_{-}$},\\
\xi_j(t+\bar t_{+}+t_{-}-\bar t_{-}-t_{+}) \quad &\text{ for $  t\geq  t_{+}+\bar t_{-}-t_{-}$}.
\end{cases}
\]
Then, by  \eqref{0-position} and \eqref{right-interval}, we have that
\[
\xi_3(0)=\eta(b_1) \quad \text{ and } \quad d(\xi_3,\xi_1)<d(\xi_2,\xi_1).
\]
This contradicts the choice of $\xi_1$ and $\xi_2$.  
Hence our claim \eqref{claim-d-1-2} holds. 

Next we show that for $k\not=l$, $(d_1(k), d_2(k))$ is not a $T$-translation of $(d_1(l), d_2(l))$. 
In fact, if 
\[
(d_1(k), d_2(k))=(d_1(l), d_2(l))+jT
\]
for some $j\in \Z\backslash \{0\}$, then both $\eta(d_1(k))$ and $\eta(d_1(l))+((k-l)J,j)$ are on $\xi_1$.  
Owing to Corollary \ref{two-integer-points},  the outward normal vector $\vec{n}$ is rational, which contradicts our assumption. 
Accordingly,  after we translate all $(d_1(k), d_2(k))$ into $(a_1, a_1+T)$,  they are disjoint.  
Therefore, \eqref{eq:sum-k} holds true.
\end{proof}

\end{document}